\documentclass{article}
\usepackage[utf8]{inputenc}

\usepackage[english]{babel}

\usepackage[letterpaper,top=2cm,bottom=2cm,left=3cm,right=3cm,marginparwidth=1.75cm]{geometry}

\usepackage{amsmath}                     
\usepackage{placeins}
\usepackage{enumitem}
\usepackage{amssymb} 
\usepackage{graphicx}
\usepackage{placeins}
\usepackage[colorlinks=true, allcolors=black]{hyperref}
\usepackage[toc,page]{appendix}
\usepackage{ulem}
\usepackage{titling}
\usepackage{xcolor}
\usepackage{xcolor}

\usepackage[affil-it]{authblk}     
\title{{Spectral Stability Correspondence between Networks and Continuous Media: Theory and Applications to Population Dynamics}}

\author{Idan Sorin}
\author{Alexander Nepomnyashchy}

\affil{Department of Mathematics, Technion–Israel Institute of Technology,\\ Haifa 32000, Israel}

\usepackage[mathlines]{lineno}

\begin{document}
\maketitle
\begin{abstract}
We investigate the stability of synchronized oscillations in coupled nonlinear systems by establishing a spectral correspondence between continuous linear shift-invariant (LSI) media and discrete networks. In this framework, Fourier modes of a continuous spatial operator and eigenmodes of a network coupling matrix are treated as spectral parameters of the same Master Stability Function.

This correspondence allows finite-wavenumber instabilities of continuous media to be translated into predictable instability windows in network coupling space. Applying the framework to zero-row-sum Metzler coupling matrices and using a competitive Lotka-Volterra model as a paradigm, we show that synchronization may exhibit reentrant behavior: it is stable for weak coupling, lost within intermediate coupling intervals, and restored at stronger coupling.

The framework also reveals a distinction between undirected and directed networks. For undirected networks, the relevant spectra are real and the resulting instability mechanism is analogous to that of standard reaction-diffusion systems with real wavenumbers. Directed networks, however, can possess complex spectra. We show that such complex spectral modes can induce quasiperiodic bifurcations of the synchronized state, leading to dynamical regimes that are inaccessible to standard real-wavenumber reflection-invariant  reaction-diffusion models.
\end{abstract}

\vspace{0.5em}
\noindent\textbf{Keywords:} {Networks, Synchronization, Master Stability Function, Instability, Lotka–Volterra model, Floquet theory, Quasiperiodic bifurcation}

\begin{quotation}
Networks of coupled nonlinear oscillators often synchronize, but the loss of
synchrony can occur in counterintuitive ways. In particular, stronger coupling
does not always make a synchronized state more stable; it may instead create
intervals of coupling strengths where synchronization is lost and then
recovered. This work explains such behavior from a spectral point of view. We
show that stability mechanisms in continuous spatial systems and in discrete
networks can be compared within a common spectral framework. For directed networks, complex spectral modes
can also change the nature of the instability, leading to quasiperiodic rather
than period-doubled desynchronization.
\end{quotation}

\section{Introduction}

The collective synchronization of coupled nonlinear oscillators
\cite{Pikovsky2001,Afraimovich1995} is a fundamental phenomenon observed
across diverse systems, ranging from neural networks and cardiac rhythms to
power grids \cite{Strogatz2003,Arenas2008}. In the context of neural dynamics,
Afraimovich, Zaks, and Rabinovich studied heteroclinic coordination in a
coupled dynamical model of sequential episodic memory initiation
\cite{Afraimovich2018}.
A central problem is to determine whether a synchronized state remains stable under perturbations transverse to the synchronization manifold. 
The Master Stability Function (MSF) formalism \cite{Pecora1998} provides a powerful framework for this problem by separating the local nodal dynamics from the spectral properties of the coupling topology \cite{Huang2009}.

While the MSF provides a stability criterion for general network topologies, it is often used as a diagnostic tool that requires case-by-case inspection of the network spectrum. 
Existing theoretical intuition is mostly limited to symmetric or undirected coupling matrices with real spectra \cite{Arenas2008,Pikovsky2016}. While
directed and non-diagonalizable networks have been incorporated into the MSF
framework \cite{NishikawaMotter2006}, the way complex spectral directions
organize stability boundaries and bifurcation mechanisms remains less
transparent. 
This analytical gap makes it difficult to predict how global parameters, such as coupling strength, determine the synchronization state. 
One might expect stronger coupling to enhance synchronization. 
However, in strongly nonlinear regimes this intuition can fail: increasing the coupling may destabilize a synchronized state, and stability may later be restored at stronger coupling values. 
Related phenomena, such as oscillation death or coupling-induced desynchronization, have been reported in nonlinear networks \cite{Ehrich2013}. 
Nevertheless, a predictive framework describing when synchronization is lost and re-established, thereby creating distinct instability windows, remains incomplete for general network topologies.

This difficulty is closely related to the historical separation between network synchronization theory and the analysis of continuous extended media. 
In continuous domains, instabilities in reaction-diffusion systems are characterized by dispersion relations, which provide analytical conditions for spatial pattern formation \cite{Cross1993,Murray2002}. 
Several works have explored connections between continuous and discrete descriptions, including the adaptation of Turing pattern analysis to discrete graphs \cite{Nakao,Wolfrum2012Turing,Asllani2014,Muolo2023} and the use of continuum limits to approximate the dynamics of large, dense networks \cite{Medvedev2014,Chiba2023}. 
However, a direct spectral correspondence that transfers instability mechanisms from continuous media to discrete network spectra is still lacking.

In this paper, we develop a spectral stability framework that bridges this gap. 
The central observation is that Fourier modes of a continuous LSI operator and transverse eigenmodes of a network coupling matrix enter the corresponding variational equations in the same spectral form. 
Although this correspondence is mathematically bidirectional, the emphasis of the present work is predictive: we use finite-wavenumber instabilities of continuous LSI media to locate synchronization and desynchronization windows in network coupling space. 
Thus, the MSF is used not only as a diagnostic tool for a prescribed network, but also as a bridge that transfers instability mechanisms from continuous dispersion relations to discrete network spectra.

We focus in particular on systems governed by zero-row-sum Metzler coupling matrices, which naturally arise in conservative transport processes and population dynamics. 
The algebraic structure of these matrices imposes strong constraints on the admissible network spectra and thereby reduces the relevant region of the MSF stability plane. 
This makes it possible to characterize directed networks with complex spectra in a systematic way.

As a representative example, we apply the framework to a competitive three-species Lotka-Volterra model. 
For real spectral modes, the finite-wavenumber instability of the corresponding continuous reaction-diffusion extension maps onto predictable desynchronization windows in the network coupling space. 
Consequently, the synchronized network state may exhibit reentrant behavior: it is stable for weak coupling, loses stability over intermediate coupling intervals, and regains stability at stronger coupling.

Complex spectral modes lead to a second, genuinely discrete-network effect. 
When the relevant network eigenvalues are complex, the instability threshold is crossed by a complex Floquet multiplier. 
The resulting bifurcation is therefore quasiperiodic rather than period-doubling, a mechanism that has no counterpart in standard reaction-diffusion systems with real wavenumbers.

The main contributions of this work are as follows:
\begin{enumerate}[label=(\roman*)]
    \item We formulate a spectral correspondence between a class of continuous LSI media and discrete network systems, showing that the Master Stability Function can be interpreted as a shared spectral stability plane for comparing continuous and discrete stability mechanisms.

    \item For zero-row-sum Metzler coupling matrices, which are relevant to conservative transport and population dynamics, we identify spectral constraints that restrict the relevant network eigenvalues to a specific sector of the left complex half-plane.

    \item Using the competitive three-species Lotka-Volterra model as a representative example, we show that finite-wavenumber instabilities of a continuous reaction-diffusion extension generate predictable synchronization and desynchronization windows in network coupling space. This leads to reentrant synchronization: increasing the coupling strength may first preserve synchronization, then destroy it over intermediate intervals, and finally restore it at stronger coupling, with multiple such transitions possible in large networks.

    \item For the same model, we demonstrate that directed networks with complex spectra produce instability mechanisms that are absent in standard reflection-invariant reaction-diffusion systems with real wavenumbers, including quasiperiodic bifurcations of the synchronized state.
\end{enumerate}

The paper is organized as follows: Section 2 establishes the general operator framework and the spectral correspondence between networks and continuous media. In Section 3, the Lotka-Volterra model is introduced and its local dynamics is analyzed. Section 4 presents the numerical results, verifying the existence of instability windows. Section 5 details the Master Stability Function analysis using the properties of Metzler matrices. Finally, concluding remarks are provided in Section 6.

\section{Spectral Correspondence: From Continuous Media to Discrete Networks}
\label{Sec2}
\subsection{Continuous Linear Shift-Invariant (LSI) Systems}
\label{sec:continuous_LSI} 

We first consider a system evolving in a continuous spatial domain $\mathbb{R}^{d}$, where the dynamics of the state vector $\mathbf{u}(\mathbf{x},t)\in\mathbb{R}^{n}$ is governed by a combination of local nonlinear reaction kinetics $\mathbf{f}(\mathbf{u})$ and a linear spatial coupling operator $\mathcal{L}$:
\begin{equation}
    \partial_{t}\mathbf{u} = \mathbf{f}(\mathbf{u}) + \mathcal{L}[\mathbf{u}].
    \label{eq:LSI_system}
\end{equation}
Here, $\mathcal{L}$ is assumed to be a Linear Shift-Invariant (LSI) operator acting on the spatial coordinates $\mathbf{x}$. According to the representation theorem for translation-invariant operators on spaces of vector-valued functions, any such operator can be uniquely represented as a convolution with a matrix-valued kernel (distribution) $G(\mathbf{x})$ \cite{Hytonen2003}:
\begin{equation}
    \mathcal{L}[\mathbf{u}](\mathbf{x},t) = (G*\mathbf{u})(\mathbf{x},t) = \int_{\mathbb{R}^{d}} G(\mathbf{x}-\mathbf{x}^{\prime})\mathbf{u}(\mathbf{x}^{\prime},t) d\mathbf{x}^{\prime}.
    \label{eq:convolution}
\end{equation}
In many ecological and biological applications, the spatial coupling typically represents transport processes (such as diffusion or nonlocal dispersal) that act on each species independently but with potentially different rates. In such cases, the operator simplifies to $\mathcal{L}[\mathbf{u}]=\mathbf{D}(g*\mathbf{u})$, where $\mathbf{D}=\text{diag}(d_{1},\dots,d_{n})$ is a diagonal matrix of diffusion coefficients, and $g(\mathbf{x})$ is a scalar kernel defining the spatial range of interaction.

Let $\mathbf{u}_{0}(t)$ be a spatially homogeneous solution of \eqref{eq:LSI_system}. To analyze its stability, we consider a perturbation of the form:
\begin{equation}
\label{perturbation}
\mathbf{u}(\mathbf{x},t)=\mathbf{u}_{0}(t)+\epsilon \mathbf{U}(t)e^{i\mathbf{k}\cdot \mathbf{x}},
\end{equation}
where $\mathbf{k}\in\mathbb{R}^{d}$ is the wave vector. The linearization of \eqref{eq:LSI_system} yields the variational equation:
\begin{equation}
    \dot{\mathbf{U}}(t) = [D\mathbf{f}(\mathbf{u}_{0}(t)) + \hat{g}(\mathbf{k})\mathbf{D}]\mathbf{U}(t),
    \label{eq:continuous_variational}
\end{equation}
where $D\mathbf{f}$ is the Jacobian of the reaction term, and $\hat{g}(\mathbf{k})$ is the Fourier transform of the kernel $g(\mathbf{x})$. The term $\hat{g}(\mathbf{k})$ represents the continuous spectrum of the spatial operator in the Fourier domain. 

\subsection{Extension to Network Systems}
\label{sec:network_extension}

The framework above naturally extends to discrete spatial domains. We consider a network of $m$ coupled nodes (patches), where the state of the $j$-th node is $\mathbf{u}_{j}\in\mathbb{R}^{n}$. The dynamics is given by:
\begin{equation}
    \frac{d\mathbf{u}_{j}}{dt} = \mathbf{f}(\mathbf{u}_{j}) + \sum_{l=1}^{m} C_{jl}\mathbf{\tilde{D}}\mathbf{u}_{l}, \quad j=1,\dots,m.
    \label{eq:network_system}
\end{equation}
Here, $C$ is an $m\times m$ coupling, or connectivity, matrix. It describes the
network topology and satisfies the zero-row-sum condition
$\sum_{l=1}^m C_{jl}=0$, which implies a zero eigenvalue with eigenvector
$(1,\ldots,1)^T$. In the context of population dynamics, this condition ensures
that uniform population distributions are not changed by migration, so that the
synchronized manifold is invariant. The matrix $\mathbf{\tilde{D}}=\text{diag}(\tilde{d}_{1},\dots,\tilde{d}_{n})$ allows for species-specific coupling strengths.

 We analyze the stability of the synchronous solution $\mathbf{u}_{1}(t)=\dots=\mathbf{u}_{m}(t)=\mathbf{u}_{0}(t)$ (where \(u_0(t)\) may be either time-dependent, for example periodic, or stationary). By defining the small perturbation $\delta u_j(t)=u_j(t)-u_0(t)$ and, when $C$ is diagonalizable, projecting it onto the eigenmodes of the coupling matrix $C$, the linearized dynamics decouples into $m$ independent equations \cite{Pecora1998, Pikovsky2016}:
\begin{equation}
    \frac{d\boldsymbol{\xi}^{s}}{dt} = [D\mathbf{f}(\mathbf{u}_{0}(t)) + \lambda^{s}\mathbf{\tilde{D}}]\boldsymbol{\xi}^{s}.
    \label{eq:network_variational}
\end{equation}

The derivation above assumes that the coupling matrix $C$ is diagonalizable. 
The extension to defective, non-diagonalizable coupling matrices can be formulated in terms of Jordan chains and is presented in Appendix~\ref{Appendix A}.

To determine stability, we consider the maximum Lyapunov exponent of Eq. \eqref{eq:network_variational} as a function of a generic complex coupling vector $\Omega=\lambda^s (\Tilde{d_1}, ...,\Tilde{d_n})^T$. This function is defined as the \textbf{Master Stability Function (MSF)} \cite{Arenas2008, Pikovsky2016}.
Equation \eqref{eq:network_variational} is the fundamental equation of the Master Stability Function (MSF) formalism. Note that \eqref{eq:network_variational} is formally equivalent to \eqref{eq:continuous_variational} under the mapping $\hat{g}(\mathbf{k})\mathbf{D} \leftrightarrow \lambda^{s}\mathbf{\tilde{D}}$. This spectral substitution shows that the MSF serves not only as a tool for network synchronization, but also as a common linear stability function for continuous LSI systems and discrete network systems. 
It is important to emphasize that this correspondence is a linear spectral correspondence. 
It applies to spatially homogeneous or synchronized base states and to their linear stability with respect to perturbations decomposed over the spectrum of the corresponding coupling operator. 
Thus, Fourier modes in the continuous system and transverse eigenmodes in the network system play analogous roles at the level of the variational equations. 
This does not imply a nonlinear equivalence between the continuous and discrete systems. 
Rather, it provides a common spectral stability plane on which the admissible spectral values of the two systems can be compared.

At the level of linear stability, choosing a network topology, through the matrix \(\mathbf{C}\), determines which points of this spectral stability plane are sampled. 
This is analogous to the way a spatial kernel \(g\) determines the Fourier symbol \(\hat{g}(\mathbf{k})\) sampled by continuous modes. 
However, the admissible spectral values are different in the two settings. 
In standard diffusive models, \(\hat{g}(\mathbf{k})\) is real, for example \(\hat{g}(\mathbf{k})=-|\mathbf{k}|^2\). 
More general continuous LSI operators, for example operators involving advection or asymmetric nonlocal kernels, may also possess complex Fourier symbols \cite{Cross1993,Kuramoto1984,Aranson2002}.
Nevertheless, their admissible spectral values remain constrained by the prescribed curve traced by \(\hat{g}(\mathbf{k})\) as the real wavenumber \(\mathbf{k}\) varies.
By contrast, directed network topologies can place transverse eigenvalues at isolated points in the complex stability plane, subject to graph-theoretic constraints. 
Thus, directed networks can sample regions of the MSF plane that are inaccessible to standard reaction-diffusion systems with real wavenumbers. 
In the present paper, we focus on applications with real \(\hat{g}(\mathbf{k})\), while future directions concerning systems with complex \(\hat{g}(\mathbf{k})\) are outlined in the conclusion.

\subsection{Spectral Constraints in Ecological Networks}
\label{Spectral Constraints in Ecological Networks}
In ecological applications where variables represent species densities, the coupling term typically models conservative flows such as migration or dispersal. In this scenario, individuals leaving a node must move to another connected node, and a node can only receive inflow from its neighbors. Consequently, the coupling matrix $C$ satisfies the properties of a zero-row-sum Metzler matrix \cite{Cvetkovic}: its off-diagonal elements are non-negative ($C_{jl} \ge 0$ for $j \neq l$, representing  inflow), and its diagonal elements are non-positive ($C_{jj} \le 0$, representing outflow). This algebraic structure is fundamental not only for population dynamics, but also appears in the modeling of resistive and capacitive electrical networks and consensus protocols \cite{Dorfler2018}.

The zero-row-sum condition implies that the diagonal term balances the off-diagonal terms:
\begin{equation}
    C_{jj} = -\sum_{l \neq j} C_{jl} \ , \ j=1,2,...,m.
\end{equation}
According to the Gershgorin circle theorem \cite{Gershgorin1931}, every eigenvalue $\lambda$ of $C$ must lie within the union of the disks defined by $|z - C_{jj}| \le |C_{jj}|$. Since $C_{jj} \le 0$, these disks are entirely contained in the left complex half-plane and pass through the origin. Therefore, we have $\text{Re}(\lambda^s) \le 0$ for all eigenvalues.

We assume that the network contains a directed spanning tree (i.e., there exists at least one root node from which all other nodes are reachable via a directed path). This assumption is not restrictive; if it does not hold, the network effectively decomposes into disjoint components or independent clusters, each of which can be analyzed as a separate subsystem. Under the zero-row-sum condition, $\lambda^1 = 0$ is a simple eigenvalue corresponding to the synchronization manifold. All transverse eigenvalues satisfy $\text{Re}(\lambda^s) < 0$, ensuring the stability of the synchronized state in the absence of reaction-induced instabilities.

Beyond the confinement to the left half-plane, the combinatorial structure of the network imposes stricter constraints on the argument of the eigenvalues. Following the analysis by Agaev and Chebotarev \cite{Agaev}, let $m$ denote the size of the matrix (number of nodes) and define the maximum off-diagonal element as $h = \max_{j \neq l} \{ C_{jl} \}$. We can represent the coupling matrix as a scaled standardized Laplacian matrix $\tilde{C}$:
\begin{equation}
    C = -h m \tilde{C}.
\end{equation}
The eigenvalues $\mu$ of the standardized matrix $\tilde{C}$ are located in the first ($Re \mu,Im \mu>0$) and fourth ($Re \mu,Im \mu<0$) quadrants. Crucially, their arguments are bounded by the network size:
\begin{equation}
    0 \le |\text{Arg}(\mu)| \le \frac{\pi}{2} - \frac{\pi}{m}.
\end{equation}
Transforming this back to the spectrum of our stability matrix $C$ (where $\lambda \propto -\mu$), this result implies that the complex eigenvalues cannot be arbitrarily close to the imaginary axis. They are confined to a sector within the left half-plane defined by the angle $\pi/m$. This constraint is particularly relevant for small networks, where the ``wedge'' of forbidden eigenvalues near the imaginary axis is significant.

Finally, since $C$ is a real matrix, its spectrum is symmetric with respect to the real axis. Since the Master Stability Function is symmetric with respect to the real axis, the stability properties for conjugate eigenvalues $\lambda^s$ and $\bar{\lambda^s}$ are identical. Consequently, it suffices to compute the Master Stability Function (MSF) only for eigenvalues in the lower-left complex half-plane (the third quadrant), subject to the angular constraints derived above.

\section{Competitive Lotka-Volterra model and its extensions}
\label{Competitive Lotka-Volterra model and its extensions}
To illustrate the spectral stability correspondence in a strongly nonlinear oscillatory setting, we consider the competitive Lotka-Volterra model, which describes the temporal evolution of $n$ biological species competing for common resources. The dynamics is governed by the system of ODEs

\begin{equation}
\frac{du_i(t)}{dt} =r_i u_i(t)\left[1-\sum_{j=1}^n \alpha_{ij}u_j (t)\right]
;\quad u_i\geq 0;\quad i=1,\ldots,n.
\label{LV}
\end{equation}

Here $u_i(t)$ is the normalized number of individuals of the $i$-th  species at time $t$, $r_i$ is the intrinsic growth rate of the $i$-th species, and $ \alpha_{ij}$ are the interaction coefficients 
that measure the extent to which the $j$-th species affects the growth rate of the $i$-th species \cite{CantrellCosner}. 
For $n=3$, the system (\ref{LV}), after a proper rescaling, takes the form:

\begin{equation}\label{LV3}
\begin{aligned}
\frac{du_1}{dt} &= r_1u_1\,(1-u_1-\alpha_1 u_2-\beta_1 u_3)\\
\frac{du_2}{dt} &= r_2u_2\,(1-u_2-\alpha_2 u_3-\beta_2 u_1)\\
\frac{du_3}{dt} &= r_3u_3\,(1-u_3-\alpha_3 u_1-\beta_3 u_2)
\end{aligned}
\end{equation}

The characteristic feature of system (\ref{LV3}) is the existence of three invariant manifolds $u_1=0$, $u_2=0$ and $u_3=0$, which bound the biologically meaningful region $u_j\geq 0$, $j=1,2,3$. For certain relations among the coefficients $\alpha_j$ and $\beta_j$, there exist heteroclinic trajectories connecting one-species saddle points in each of the manifolds, which together form a {\it robust} heteroclinic cycle \cite{Coste1979, Chi,MayLeonard}.  In \cite{Coste1979}, it has been shown that in the case of unequal growth rates $r_i$, this system can have three possible types of attractors: a coexistence state of equilibrium (none of the variables equals zero), a limit cycle, and a heteroclinic cycle. Note that in the special case 
$r_1=r_2=r_3=r$,  which is called the May-Leonard system~\cite{MayLeonard, Sorin2025}, the dynamics are different: for that system, it was proven \cite{Chi} that for parameters that satisfy the relation 
$$(1-\beta_1)(1-\beta_2)(1-\beta_3)=(\alpha_1-1)(\alpha_2-1)(\alpha_3-1),\quad 0<\beta_i<1<\alpha_i,$$ 
there exists an attracting two-dimensional invariant manifold with a continuum of periodic solutions on it.

To simplify the notation, we denote $u\equiv u_1$, $v\equiv u_2$ and $w\equiv u_3$. As a basic example, we consider in detail the following two-parameter system,
which also appears as the dissipative component in related work on seasonally
forced cyclic-competition dynamics \cite{Sorin2026}: 
\begin{subequations}
\label{specific}
\begin{align}
u_t &= u(\gamma - u - \alpha v)  ,  \\
v_t &= v(1 - v - \alpha w)  , \ \ \alpha,\gamma >0, \\
w_t &= w(1 - w - \alpha u). 
\end{align}
\end{subequations}
Here and below, a subscript $t$ or $x$ means the derivative with respect to the corresponding variable.

The basins of attraction of the system in the parameter region $0<\gamma<1$ are indicated in Figure \ref{Parameters space}. Evaluation of the exact expressions for these regions of attraction can be found in Appendix \ref{Appendix B}.
While model \eqref{specific} might look simpler than the general model \eqref{LV3}, it was shown in \cite{Coste1979} that similar dynamics also occurs in the general competitive 3-species Lotka-Volterra system, meaning that system \eqref{LV3} can also have only three possible attractors: a coexistence point, a limit cycle and a heteroclinic cycle. Thus, system \eqref{specific} provides a minimal representative setting that retains the
dynamics required for the stability mechanism studied below. We denote the periodic solution of \eqref{specific}, with period $T$ by:
\begin{equation}
\label{periodic}
(u_0(t),v_0(t),w_0(t)).
\end{equation}

The behavior near heteroclinic cycles has been analyzed in more detail in systems with four \cite{Groothuizen2022} and five \cite{Postlethwaite2022} competitive species.

\begin{figure}[htbp]
\centering
\includegraphics[width=0.7\linewidth]{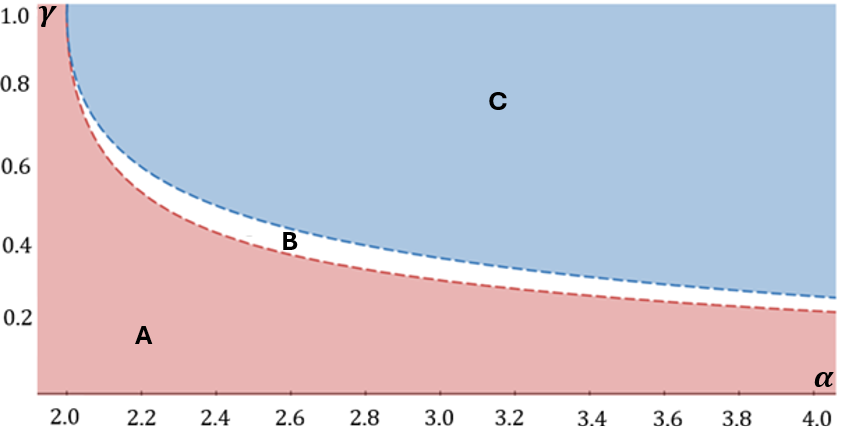}
\caption{\label{Parameters space} The parameter space $(\alpha,\gamma)$ of system \eqref{specific} (computed analytically). In region (A) the coexistence equilibrium is stable, in region (B) there is a stable} limit cycle, and in region (C) the system is attracted to a heteroclinic cycle.  
\end{figure}

In the following, we consider two extensions of the local system that describe species migration or dispersal.

The continuous reaction-diffusion extension of (\ref{specific}) is:
\begin{subequations}
\label{RD specific}
\begin{align}
u_t &= u(\gamma - u - \alpha v)+d_uu_{xx}  ,  \\
v_t &= v(1 - v - \alpha w)+d_vv_{xx}  , \\
w_t &= w(1 - w - \alpha u)+d_ww_{xx};\;
-\infty<x<\infty.
\end{align}
\end{subequations}
Here $u$, $v$ and $w$ are functions of the temporal coordinate $t$ and the spatial coordinate $x$. The system (\ref{RD specific}) has a solution that does not depend on $x$ and coincides with the solution of the ODE system discussed above:
\begin{equation}
u(x,t)=u_0(t),\;v(x,t)=v_0(t),\;w(x,t)=w_0(t),
\label{uniform}
\end{equation}
which corresponds to synchronous oscillations throughout the space.

We  shall also consider another kind of extension to (\ref{specific}), namely the network:

\begin{subequations}
\label{Network specific}
\begin{align}
\dot{u_j} &= u_j(\gamma - u_j- \alpha v_j)+D_u\sum_{l=1}^m C_{jl}u_l  ,  \\
\dot{v_j} &= v_j(1 - v_j - \alpha w_j)+D_v\sum_{l=1}^m C_{jl}v_l, \\
\dot{w_j} &= w_j(1 - w_j - \alpha u_j)+D_w\sum_{l=1}^m C_{jl}w_l. 
\end{align}
\end{subequations}
$j=1,...,m$.
 All subsystems that make up the network are assumed to be identical, and the coupling matrix is the same for all species.

Biologically, this can be understood as $m$ different habitats, with the same ecological system in each of them. In each habitat, the same species live with the same dynamics between them, but individuals of each species can migrate between habitats. 
As in Section \ref{sec:network_extension}, we assume the existence of a synchronized periodic solution with
\begin{equation}
u_j(t)=u_0(t),\;v_j(t)=v_0(t),\;w_j(t)=w_0(t)
\label{synchronized}
\end{equation}
and impose the condition $\sum_{l=1}^m C_{jl}=0$ for every $1\leq j\leq m$,
which means that the sum of every row in the connectivity matrix is zero.

\section{Instabilities of synchronous oscillations}
\label{Sec4}

The purpose of this section is to demonstrate that the spectral correspondence developed in Section~\ref{Sec2} is predictive. 
As mentioned above, the periodic solution~(\ref{periodic}) is stable in the local system~(\ref{specific}), whereas its spatial extension~(\ref{uniform}) and network extension~(\ref{synchronized}) may become unstable with respect to desynchronizing perturbations. 
Rather than performing an exhaustive parameter sweep, we focus on a representative and non-trivial mechanism: desynchronization induced by increasing coupling strength. 
In the continuous extension, this mechanism appears as a finite-wavenumber spatial instability interacting with a period-doubling bifurcation \cite{Sorin2025}. 
Through the spectral correspondence, once the finite-wavenumber instability interval of the continuous problem is known, the corresponding desynchronization intervals of the network can be obtained directly from the eigenvalues of $C$. 
For directed networks with complex spectra, the same framework also shows how this underlying instability can be transformed into quasiperiodic bifurcations, thereby producing dynamical regimes that are inaccessible to standard continuous reaction-diffusion systems with real wavenumbers.

\subsection{Period-doubling instability for finite wavenumbers in the continuous extension}

As in (\ref{perturbation}), we add a small disturbance to the time-periodic solution (\ref{periodic}):
\begin{equation}  
(u,v,w)=(u_0(t)+\epsilon U(t)e^{ikx},v_0(t)+\epsilon V(t)e^{ikx},w_0(t)+\epsilon W(t)e^{ikx}) , \ 0<\epsilon\ll 1,
\label{perturbation specific}
\end{equation}
so that the linearized system is: 
\begin{equation}
\begin{pmatrix}
\dot{U} \\
\dot{V} \\
\dot{W}
\end{pmatrix}
= - \begin{pmatrix}
d_u k^2 - \gamma + 2u_0(t)+\alpha v_0(t) & \alpha u_0(t) & 0 \\
0 & d_v k^2 - 1 + 2v_0(t)+\alpha w_0(t) & \alpha v_0(t) \\
\alpha w_0(t) & 0 & d_w k^2 - 1 + 2w_0(t)+\alpha u_0(t)
\end{pmatrix}
\begin{pmatrix}
U \\
V \\
W
\end{pmatrix}.
\label{linear specific}
\end{equation}

We aim to analyze the stability of the spatially homogeneous solution \eqref{uniform} with respect to the perturbation \eqref{perturbation specific}.
In our numerical examples, we set the     diffusion coefficients $d_v = d_w = 0$ to reduce the dimensionality of the parameter space. We note that this choice is made for simplicity of presentation. We also set $\gamma=0.5$, so $\alpha$ is the bifurcation parameter. We choose
$\alpha$ in the region (B) of Figure \ref{Parameters space}, $\alpha_-<\alpha<\alpha_+$, where $\alpha_-\approx 2.24$, $\alpha_+\approx 2.38$. 

By integrating (\ref{linear specific}) within the interval $0\leq t\leq T$ for three linearly independent initial conditions, we obtain the monodromy matrix ${\bf M}(k^2)$ that relates the vectors $(U(T),V(T),W(T))^T$ and $(U(0),V(0),W(0))^T$:
\begin{equation}
\begin{pmatrix}
U(T) \\
V(T) \\
W(T)
\end{pmatrix}
= {\bf M}(k^2)
\begin{pmatrix}
U(0) \\
V(0) \\
W(0)
\end{pmatrix}.
\label{monodromy}
\end{equation}
The eigenvalues $\mu_i$, $i=1,2,3$, of the monodromy matrix (multipliers) determine the stability of the solution (\ref{periodic}).

\begin{figure}[!t]
\centering
\includegraphics[width=0.9\linewidth]{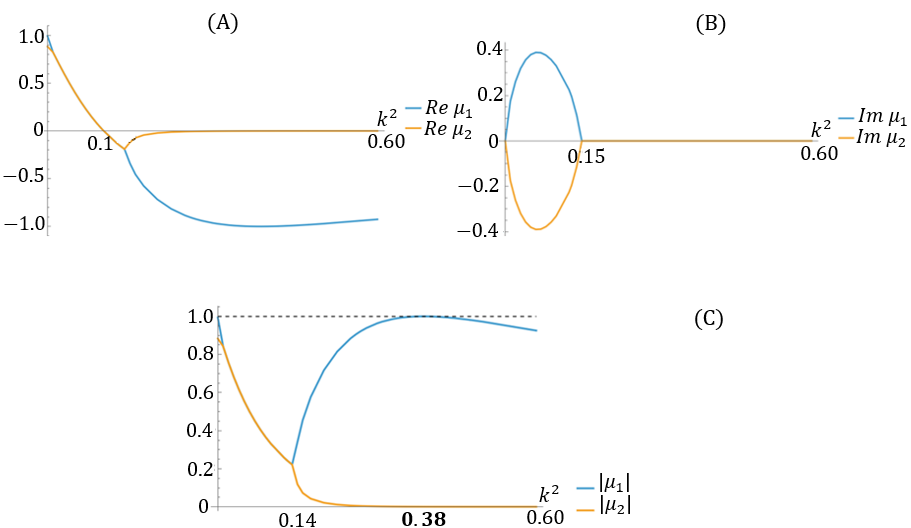}
\caption{\label{PD1} Period doubling of the solution \eqref{perturbation specific} at a finite wavenumber for $d_u=1,d_v=d_w=0$. The period doubling occurs at $k^2=k_*^2 \approx 0.38$. In (A),(B),(C) the real parts, imaginary parts, and moduli of the two leading Floquet
multipliers are shown respectively at $\alpha=2.3427$.}
\end{figure}

\begin{figure}[!t]
\centering
\includegraphics[width=0.9\linewidth]{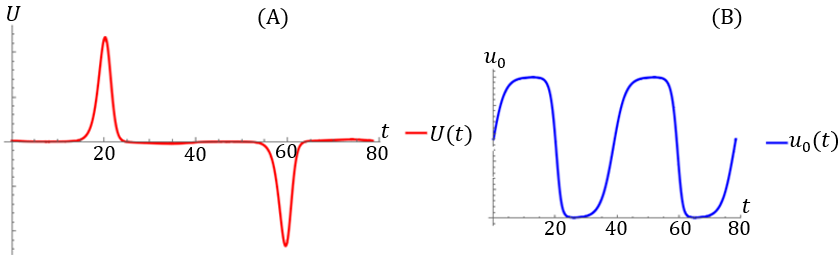}
\caption{\label{PD2} Period doubling of the solution \eqref{perturbation specific} at a finite wavenumber, computed numerically for $\gamma=0.5$},$d_u=1, d_v=d_w=0$. (A) The disturbance $u$-component, denoted as $U(t)$, at the critical wavenumber $k=k_*\approx 0.62$ and $\alpha=2.3427$. (B) The unperturbed base solution $u_0(t)$. Notice that the period of $U(t)$ in (A) is exactly twice the period of the base solution in (B).
\end{figure}

All numerical integrations and computations of the Floquet multipliers were performed using Wolfram Mathematica, utilizing specialized solvers designed for stiff differential equations to ensure numerical stability and accuracy.
The results of the multiplier calculation for $\alpha=2.3427,$ $\gamma=0.5,d_u=1$ are shown in Figure \ref{PD1}. One of the multipliers remains strongly contracting and is therefore not shown in the plots. For the other two multipliers, Re$\mu_i$, Im$\mu_i$, and $|\mu_i|$ are shown as functions of $k^2$.
At $k^2=0$, one of the multipliers is equal to 1 (it corresponds to the time shift) 
and another is positive and smaller than 1, which confirms the stability of solution  (\ref{periodic}) in the framework of the system (\ref{specific}). At small $k^2$, both multipliers are positive and smaller than 1; This means that the synchronous oscillations are stable with respect to spatially periodic disturbances. There exists a critical value $\alpha=\alpha_{*}$ , such that at $\alpha<\alpha_*$, $|\mu(k^2)|<1$ for all multipliers for any $k^2\neq 0$, that is, the spatially uniform time-periodic solution is stable. However, at $\alpha>\alpha_*$, $\mu(k^2)<-1$ in a certain interval $k_{min}^2<k^2<k_{max}^2$. That means that within that interval the multiplier $\mu(k^2)$ crosses the value -1, creating the instability with respect to oscillations with period $2T$ and periodic in space with period $2\pi/k$. We find the critical value $\alpha=\alpha_*\approx 2.34$, and for that value $\mu(k^2)=-1$ at $k=k_*\approx 0.62$. Thus, the temporal period of the critical disturbance is twice as large as that of the basic oscillation (see Figure \ref{PD2}).
One can expect that the instability described above leads to the appearance of a spatially periodic solution with a double temporal period \cite{Argentina}. This expectation is supported by the network simulations discussed below.

The instability described above is similar to that found in \cite{Sorin2025} for periodic solutions of the May-Leonard system close to the heteroclinic   cycle. Another example of this instability is given in \cite{Argentina}.

Note that as $\alpha$ increases further, the time-periodic oscillations are replaced by quasiperiodic oscillations due to a secondary instability. The analysis of secondary bifurcations is beyond the scope of this paper.

\subsection{Network extension and desynchronization windows}

Consider now the network (\ref{Network specific}). Linearizing equations (\ref{Network specific}) around the periodic solution (\ref{synchronized}) corresponding to synchronous oscillations, we obtain for disturbances $(U_j(t),V_j(t),W_j(t))$, $j=1,\ldots,m$:
\begin{subequations}
\label{Network linear}
\begin{align}
\dot{U_j} &= (\gamma - 2u_0- \alpha v_0)U_j-\alpha u_0V_j+D_u\sum_{l=1}^m C_{jl}U_l  ,  \\
\dot{V_j} &= (1 - 2v_0 - \alpha w_0)V_j-\alpha v_0W_j+D_v\sum_{l=1}^m C_{jl}V_l, \\
\dot{W_j} &= (1 - 2w_0 - \alpha u_0)W_j-\alpha w_0U_j+D_w\sum_{l=1}^m C_{jl}W_l,
\end{align}
\end{subequations}
$j=1,...,m$. 

As explained in Section 2.3, the expansion of vectors $(U_j,V_j,W_j)$ over the eigenvectors of the coupling matrix $C$ splits the system (\ref{Network linear}) of order $3m$ into $m$ decoupled systems equivalent to (\ref{linear specific}) with $-d_uk^2$, $-d_vk^2$, $-d_wk^2$ replaced with $D_u\lambda^s$, $D_v\lambda^s$, $D_w\lambda^s$, where $\{\lambda^s\}$, $s=1,\ldots,m$, are the eigenvalues of matrix $C$.

Before studying the general case, consider several particular cases.

\subsubsection{Two coupled subsystems}
\label{Two coupled subsystems}
First, let us consider the simplest ``network" that consists of two coupled subsystems,
\begin{equation}
\dot{u}_1=u_1(\gamma-u_1-\alpha v_1)+D(u_2-u_1),\;\dot{v}_1=v_1(1-v_1-\alpha
w_1),\;\dot{w}_1=w_1(1-w_1-\alpha u_1),
\label{first_nonlinear}
\end{equation}  
\begin{equation}
\dot{u}_2=u_2(\gamma-u_2-\alpha v_2)+D(u_1-u_2),\;\dot{v}_2=v_2(1-v_2-\alpha
w_2),\;\dot{w}_2=w_2(1-w_2-\alpha u_2),
\label{second_nonlinear}
\end{equation}  
which has the synchronized basic solution
$$u_1=u_2=u_0(t),\;v_1=v_2=v_0(t),\;w_1=w_2=w_0(t).$$
Its stability is determined by the system
$$\dot{U}_1=U_1(\gamma-2u_0-\alpha v_0-D)-\alpha
u_0V_1+DU_2,\;\dot{V}_1=V_1(1-2v_0-\alpha
w_0)-\alpha v_0W_1,\;$$
\begin{equation}
\dot{W}_1=W_1(1-2w_0-\alpha u_0)-\alpha w_0U_1,
\label{first_linear}
\end{equation}
$$\dot{U}_2=U_2(\gamma-2u_0-\alpha v_0-D)-\alpha
u_0V_2+DU_1,\;\dot{V}_2=V_2(1-2v_0-\alpha
w_0)-\alpha v_0W_2,$$
\begin{equation}
\dot{W}_2=W_2(1-2w_0-\alpha u_0)-\alpha w_0U_2,
\label{second_linear}   
\end{equation}  
Define
$$U_{\pm}=U_1\pm U_2,\;V_{\pm}=V_1\pm V_2,\;W_{\pm}=W_1\pm W_2.$$
Then
$$\dot{U}_+=U_+(\gamma-2u_0-\alpha v_0)-\alpha
u_0V_+,\;\dot{V}_+=V_+(1-2v_0-\alpha
w_0)-\alpha v_0W_+,$$
\begin{equation}
\dot{W}_+=W_+(1-2w_0-\alpha u_0)-\alpha w_0U_+, 
\label{synchronous_linear}
\end{equation}
$$\dot{U}_-=U_-(\gamma-2u_0-\alpha v_0-2D)-\alpha
u_0V_-,\;\dot{V}_-=V_-(1-2v_0-\alpha
w_0)-\alpha v_0W_-,$$
\begin{equation}
\dot{W}_-=W_-(1-2w_0-\alpha u_0)-\alpha w_0U_-.
\label{asynchronous_linear}
\end{equation}
Equations (\ref{synchronous_linear}) and (\ref{asynchronous_linear}) are each equivalent to (\ref{linear specific}) with $k^2=0$ and
$d_uk^2=2D$, respectively. Therefore, there exist two branches of multipliers: one with constant multipliers,
$$\mu_j^+(D)=\mu_j(0),\;j=1,2,3,$$
equal to the multipliers of problem (\ref{linear specific}) at $k^2=0$ and a branch 
$$\mu_j^-(D)=\mu_j(2D/d_u),\;j=1,2,3,$$
with multipliers equal to the multipliers of the problem (\ref{linear specific}) at 
\begin{equation}
k^2=2D/d_u.
\label{correspondence}
\end{equation}
The plot of $\mu_j^-(D)$ for the desynchronization mode is equivalent to that shown in Figure \ref{PD1} with replacement (\ref{correspondence}).

At $\alpha>\alpha_*$, $\alpha_*\approx 2.34$, the instability of synchronous oscillations takes place in the interval of the values of the coupling constant $D_{-}<D<D_+$, where $D_{\pm}=d_uk_{\pm}^2(\alpha)/2$. While an increase in the coupling parameter is conventionally expected to promote synchronization \cite{Kuramoto1975, Strogatz2015}, we observe a counterintuitive non-monotonic dependence. Specifically, synchronous oscillations remain stable for sufficiently weak coupling, $0 < D < D_-$, but undergo desynchronization as the coupling constant $D$ is increased beyond $D_-$. Stability is only restored at much higher coupling strengths, $D > D_+$. Although desynchronization induced by strengthened coupling has been reported in specific discrete neural networks \cite{Ehrich2013}, our results establish this phenomenon within the context of ecological networks. Furthermore, as we will demonstrate in the subsequent section, this loss of stability is not a singular event; rather, the system can traverse multiple distinct windows of desynchronization whose precise locations and widths can be analytically predicted.

\begin{figure}[!t]
\centering
\includegraphics[width=0.9\linewidth]{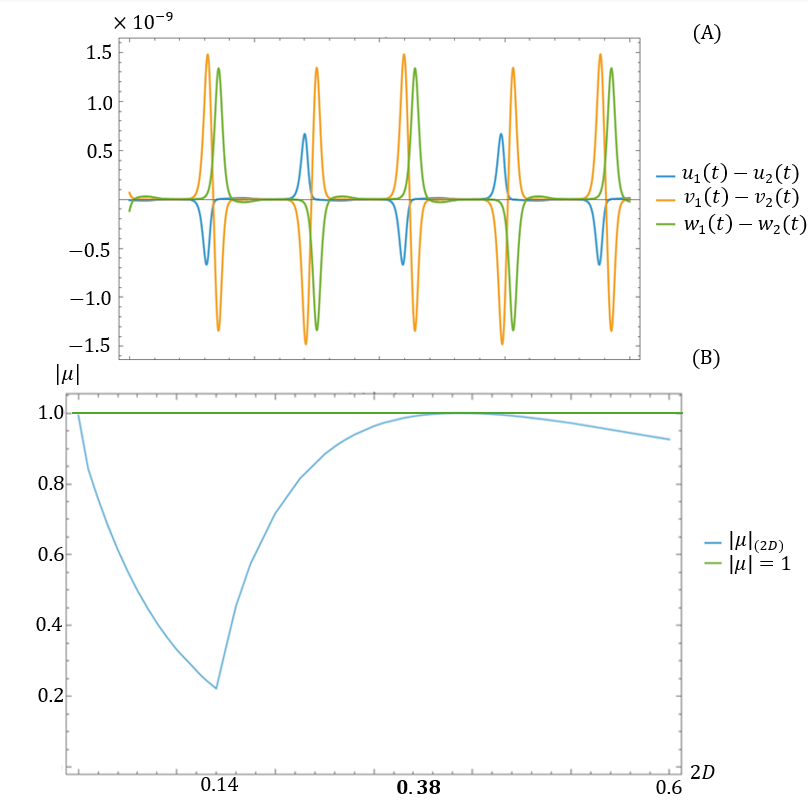}
\caption{\label{Two_systems.png} Period-doubling bifurcation in the network of two coupled subsystems \eqref{first_nonlinear} and \eqref{second_nonlinear}, computed numerically. (A) The differences between the solutions of the two subsystems as a function of time. (B) The leading Floquet multiplier $\mu$ as a function of the doubled coupling strength $2D$. Parameters are $\alpha=2.3427$ and $\gamma=0.5$, as in Fig. \ref{PD1}. The period-doubling desynchronization occurs at $D=k_*^2/2\approx 0.1922$, where the multiplier reaches $\mu=-1$.}
\end{figure}

 The period-doubling bifurcation in the two-site system is shown in Figure \ref{Two_systems.png}. Figure \ref{Two_systems.png}B shows the dependence $|\mu(D)|$ for the multiplier with the highest modulus value. Figure \ref{Two_systems.png}A demonstrates the oscillations with the double period obtained by numerical simulation of the system (\ref{first_nonlinear}), (\ref{second_nonlinear}).

\begin{figure}[!t]
\centering
\includegraphics[width=0.8\linewidth]{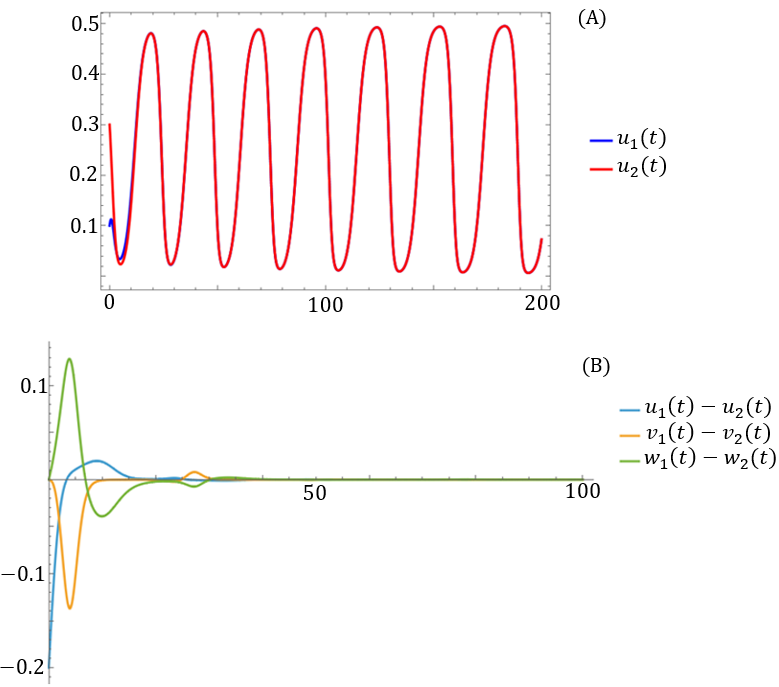}
\caption{\label{Syn.png} Synchronization between the two coupled subsystems \eqref{first_nonlinear} and \eqref{second_nonlinear} with differing initial data, computed numerically. (A) The trajectories $u_1(t)$ and $u_2(t)$ as a function of time. (B) The synchronization errors $u_1-u_2$, $v_1-v_2$, and $w_1-w_2$ as a function of time. Parameters: $\alpha=2.3427$, $\gamma=0.5$, and $D=0.15$. The initial conditions are $(u_1(0),v_1(0),w_1(0))=(0.1,0.15,0.05)$ and $(u_2(0),v_2(0),w_2(0))=(0.3,0.15,0.05)$. Notice that despite the different initial states, the errors quickly decay to zero, indicating synchronization.}
\end{figure}

In the regions of linear stability of  synchronous oscillations, direct numerical simulations of the full nonlinear system show that the oscillations synchronize in both subsystems for different initial conditions. Figure \ref{Syn.png} shows that solutions of system (\ref{first_nonlinear}), (\ref{second_nonlinear}) with different initial conditions converge to the synchronized periodic state with $u_1(t)=u_2(t),v_1(t)=v_2(t),w_1(t)=w_2(t)$.

\subsubsection{Real network spectra and desynchronization windows}
Now, consider a network of more than two coupled subsystems, but assuming that all the eigenvalues of the matrix $C$ are real and negative. This is always the case when $C$ is symmetric and Metzler. It was shown in Section \ref{Spectral Constraints in Ecological Networks} that one of the eigenvalues is zero and the others are strictly negative. 

Assume that for a set of parameters $(d_u,d_v,d_w,k)$, the stability analysis for the PDE model \eqref{RD specific} predicts the instability in a certain interval $k^2_{min}<k^2<k^2_{max}$. Let us consider system \eqref{Network linear} and start changing $D_u$, $D_v$ and $D_w$ proportionally keeping constant the ratios $D_v/D_u=d_v/d_u$, $D_w/D_u=d_w/d_u$. For each eigenvalue $\lambda^s <0$, $s=1,\ldots,m$, of the matrix $C$, three multipliers of system \eqref{Network linear} for a given value of $D_u$ coincide with three multipliers of system \eqref{linear specific} at the value of $k$ satisfying the relation $D_u\lambda^s=-d_uk^2$. Thus, for each $\lambda^s$, there will be instability in the interval $-d_uk_{min}^2/\lambda^s<D_u<-d_uk_{max}^2/\lambda^s$.

\begin{figure}[htbp]
\centering
\includegraphics[width=0.98\linewidth]{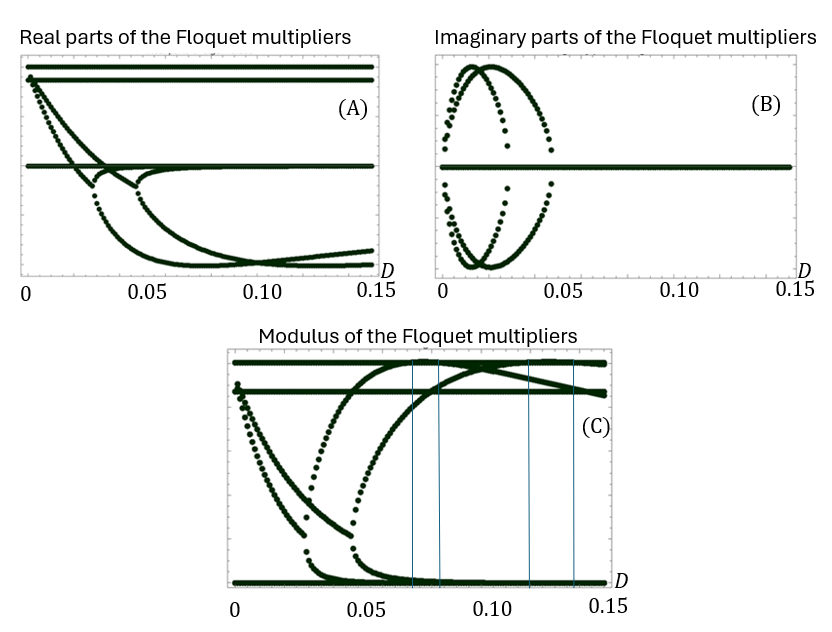}
\caption{\label{two_intervals.png} 
Floquet multipliers for the linear system~\eqref{linearlambda}, where the coupling matrix \(C\) has eigenvalues \(\{-5,-3,0\}\), computed numerically. 
Panels (A), (B), and (C) show the real parts, imaginary parts, and moduli of the multipliers, respectively, as functions of the coupling strength \(D\). 
Parameters are \(\alpha=2.3435\), \(\gamma=0.5\), \(D_u=D\), and \(D_v=D_w=0\). 
The key information is shown in panel (C): the vertical lines mark the boundaries of the two intervals in which at least one Floquet multiplier satisfies \(|\mu|>1\). 
These intervals, \(0.068<D<0.088\) and \(0.113<D<0.147\), correspond to two distinct desynchronization windows generated by the two nonzero transverse eigenvalues of \(C\). 
Thus, as predicted by the spectral correspondence, the synchronized state is lost and then recovered twice as the coupling strength is increased. 
All \(9\) multipliers of the full linearized network system are plotted, although some of them are degenerate.}
\end{figure}

When two eigenvalues \(|\lambda^{s_1}|<|\lambda^{s_2}|\) are sufficiently close, namely when
\[
\left|\frac{\lambda^{s_2}}{\lambda^{s_1}}\right| < \frac{k^2_{\max}}{k^2_{\min}},
\]
the corresponding instability intervals overlap and merge into a single larger desynchronization window. 
If, on the other hand,
\[
\left|\frac{\lambda^{s_2}}{\lambda^{s_1}}\right| > \frac{k^2_{\max}}{k^2_{\min}},
\]
the two eigenvalues generate distinct instability windows. 
Thus, as the coupling strength increases, the network may lose and regain synchronization more than once. 
The number and width of these windows can be controlled through the spectrum of the coupling matrix \(C\). Consider a network with non-positive eigenvalues ordered by magnitude as $0=|\lambda^1|<|\lambda^2|<\dots <|\lambda^m|$. If the spectrum satisfies the separation condition $|\lambda^{j+1}/\lambda^{j}| > k^2_{max}/k^2_{min}$ for $2 \leq j\leq m-1$, the system will exhibit exactly $m-1$ distinct desynchronization windows. These instability regions correspond to the intervals $|d_u k^2_{min}/\lambda^j| < D_u < |d_u k^2_{max}/\lambda^j|$, while intermediate synchronization windows emerge in the ranges $|d_u k^2_{max}/\lambda^j| < D_u < |d_u k^2_{min}/\lambda^{j+1}|$. An example of a network with two desynchronization windows is shown in Figure \ref{two_intervals.png} using the coupling matrix
$$C=
\begin{bmatrix}
-3 & 1 & 2 \\
1 & -2 & 1 \\
2 & 1 & -3
\end{bmatrix}.$$\\

\subsubsection{Complex network spectra and quasiperiodic desynchronization}
\label{Networks with complex eigenvalues}

We now consider the case where the coupling matrix $C$ possesses a pair of complex conjugate eigenvalues $\lambda^s$ and $\bar{\lambda}^s$. This scenario arises naturally in directed networks with non-symmetric coupling. Note that the stability results obtained within the standard PDE extension are insufficient here, as they rely strictly on real wavenumbers $k$. However, as shown in Section \ref{Spectral Constraints in Ecological Networks}, for zero-row-sum Metzler matrices, the spectrum is confined to the left half of the complex plane; thus, it suffices to analyze $\lambda^s$ solely in the second or third quadrant.

We assume $D_u \neq 0$ and seek the Floquet multipliers for the decoupled mode equations:
\begin{subequations}
\label{linearlambda}
\begin{align}
    \dot{U}^s &= (\gamma-2u_0-\alpha v_0)U^s-\alpha u_0V^s+D_u\lambda^sU^s, \label{Du} \\
    \dot{V}^s &= (1-2v_0-\alpha w_0)V^s-\alpha v_0W^s+D_v\lambda^sV^s, \label{Dv} \\
    \dot{W}^s &= (1-2w_0-\alpha u_0)W^s-\alpha w_0U^s+D_w\lambda^sW^s, \label{Dw}
\end{align}
\end{subequations}
as functions of the complex parameter $\Omega \equiv D_u\lambda^s$, while maintaining constant ratios $D_v/D_u$ and $D_w/D_u$.

First, consider the regime where $|\Omega|$ is small. According to Floquet theory, the monodromy matrix takes the form $M=\exp(\Lambda T)$, where $T$ is the period of the unperturbed limit cycle $(u_0(t),v_0(t),w_0(t))$. At $\Omega=0$, the matrices $M$ and $\Lambda$ possess distinct real eigenvalues $\mu_i(0)$ and $\Lambda_i(0)$, satisfying $\mu_i=\exp(\Lambda_i T)$. The largest multiplier corresponds to the time-translation mode, with $\mu_1(0)=1$ and $\Lambda_1(0)=0$.
From the PDE analysis (where diffusion corresponds to real negative $\Omega$), it is known that for small real negative $\Omega$, the leading multiplier satisfies $\mu_1(\Omega)<1$, implying a negative growth rate $\Lambda_1(\Omega)<0$. Consequently, $\mu_1$ admits a series expansion: $\mu_1(\Omega) = 1 + a_1\Omega + O(\Omega^2)$, with $a_1 > 0$. By continuity, for complex $\Omega$ with a sufficiently small imaginary part and $\text{Re}(\Omega) < 0$, we have:
\begin{equation}
    \text{Re}(\Lambda_1(\Omega)) < 0.
\end{equation}
Therefore, the modulus of the multiplier satisfies:
\begin{equation}
    |\mu_1(\Omega)| = |\exp(\Lambda_1(\Omega))| = \exp(\text{Re}(\Lambda_1(\Omega))) < 1.
\end{equation}
This confirms that no desynchronizing instability occurs for sufficiently small $\Omega$ in the left half-plane. In other words, if the continuous PDE limit predicts stability against long-wave perturbations, the network remains stable at sufficiently weak coupling.

\begin{figure}[htbp]
    \centering
    \includegraphics[width=0.7\linewidth]{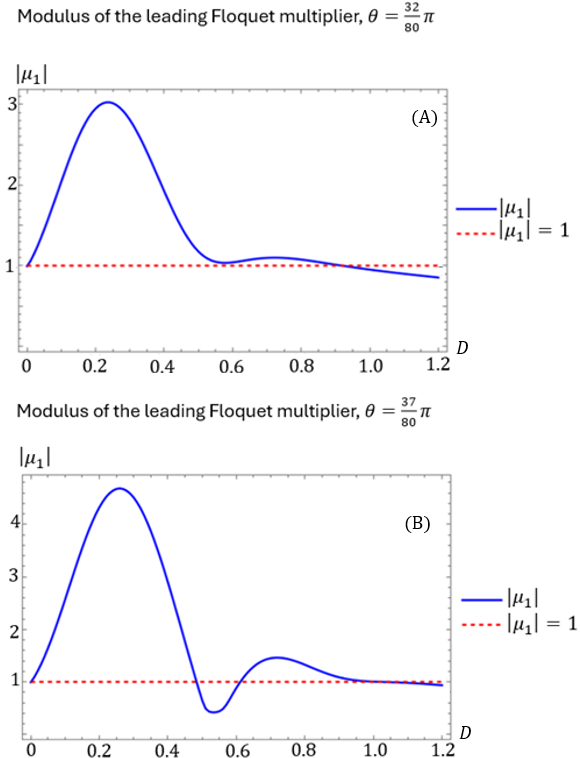}
    \caption{\label{complex_eigenvalues} The modulus (absolute value) of the leading Floquet multiplier $|\mu_1|$ for the system \eqref{linearlambda} with complex network eigenvalues $\lambda^s=-e^{i\theta}$, computed numerically. (A) For $\theta =\frac{32}{80}\pi$, the local minimum of $|\mu_1|$ lies above the threshold $|\mu_1|=1$ (dashed line). Consequently, there is a single, continuous instability window for $\epsilon < D < 0.92$ (where $\epsilon > 0$ is too small to be visible in the figure). (B) For $\theta =\frac{37}{80}\pi$, the local minimum drops below $|\mu_1|=1$, splitting the instability region into two distinct intervals: $\epsilon^* < D < 0.49$ and $0.61 < D < 1.04$.}
    \label{BigTheta}
\end{figure}

We now focus on the numerical results for the parameters $\alpha=2.3427$ (close to the bifurcation at $\alpha_*$) and $D_v=D_w=0$, denoting $D_u=D$. We parameterize the eigenvalue as $\lambda^s = -|\lambda^s|e^{i\theta}$.
Similarly to the case of real eigenvalues, the synchronized solution is stable at both weak and strong coupling limits.
The behavior of the leading multiplier $|\mu_1|$ as a function of $D$ depends on the angle $\theta$:
\begin{itemize}
    \item For small $\theta$, $|\mu_1(D)|$ exhibits a single maximum and minimum, similar to the real case, creating a single interval of instability.
    \item As $\theta$ increases beyond a critical value, the curve $|\mu_1(D)|$ deforms: the inflection point splits, creating complex stability landscapes.
    \item Eventually, for larger $\theta$, the local minimum may drop below 1 while the maximum remains above 1, leading to two disjoint instability windows.
\end{itemize}
These scenarios are illustrated in Figure \ref{BigTheta}.

\begin{figure}[htbp] 
    \centering
    \includegraphics[width=1.0\linewidth]{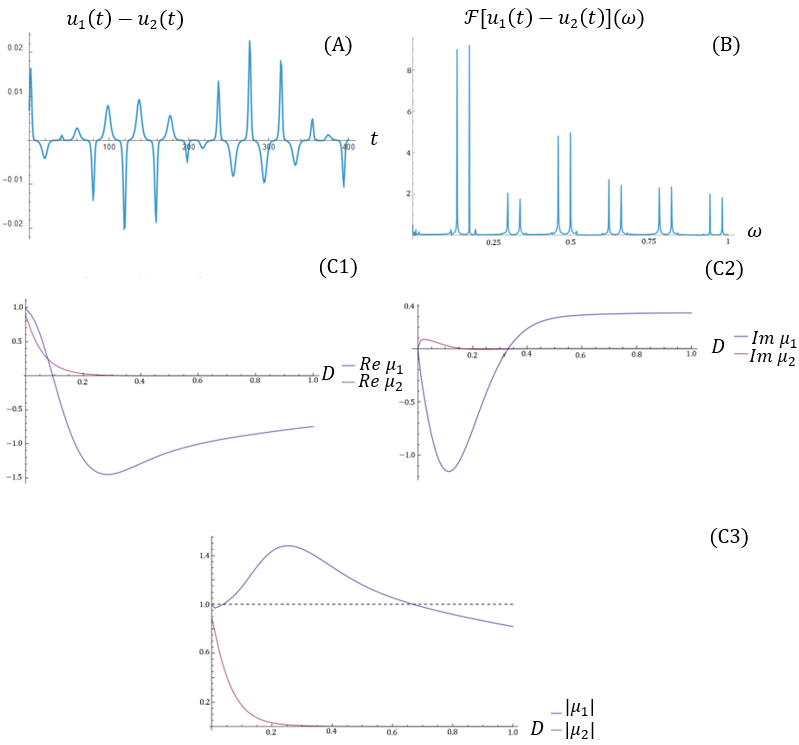}
    \caption{\label{quasi_periodic} Quasiperiodic solution of the system \eqref{Network specific}, computed numerically. (A) The difference $u_1-u_2$ between the two subsystems as a function of time. (B) The Fourier transform of the difference, showing two basic frequencies $\omega_1$ and $\omega_2$. Panels (C1), (C2), and (C3) show the real parts, imaginary parts, and moduli of the two leading Floquet multipliers $\mu_1,\mu_2$ as a function of the coupling strength $D$, respectively. The network is defined by the connectivity matrix $C$ (shown below), which has the eigenvalues $\{-\sqrt{2},\,0,\,-e^{\pm i\pi/4}\}$. Parameters: $D_u=D$ and $D_v=D_w=0$, with $D=0.04$ in panels (A) and (B).}
\[
C=\frac{1}{\sqrt{2}}
\begin{pmatrix}
-1 &  1 &  0 &  0 \\
 0 & -1 &  1 &  0 \\
 0 &  0 & -1 &  1 \\
 1 &  0 &  0 & -1
\end{pmatrix}
\]

    \label{QP}
\end{figure}

Crucially, for $\theta \neq 0$, the critical multiplier $\mu_1$ is generally complex at the instability threshold ($|\mu_1|=1$). This implies that instead of a simple period-doubling bifurcation (characteristic of real spectra), the system undergoes a quasiperiodic (Neimark-Sacker) bifurcation \cite{Sacker2009Chapter}.
Figure \ref{QP} demonstrates the resulting quasiperiodic oscillations of the synchronization error $u_1(t)-u_2(t)$ just above the instability threshold. This dynamical regime has no analog in reaction-diffusion systems or undirected networks with symmetric coupling, since this is a network-induced bifurcation generated by the complex direction in the MSF plane.

The quasiperiodic nature is evident from the spectral decomposition. Since the imaginary part of the critical multiplier is non-zero, the solution takes the form:
\begin{equation}
    u_j(t) = u_0(t) + \epsilon \sum_{n,m \in \mathbb{Z}\setminus\{0\}} A^{(j)}_{n,m} e^{i(n\omega_1 + m\omega_2)t}
\end{equation}
(and similarly for $v_j, w_j$). The Fourier transform $\mathcal{F}[u_1(t)-u_2(t)](\omega)$ reveals two incommensurate fundamental frequencies:
\begin{equation}
    \omega_1 = \frac{2\pi}{T}, \quad \omega_2 = \frac{\text{Arg}(\mu_1)}{T}. 
\end{equation}
Numerically, we find $\omega_1 \approx 0.16$ and $\omega_2 \approx 0.019$. Since zeroth-order terms cancel out in the difference $u_1-u_2$, the spectrum is dominated by peaks at combinational frequencies $n\omega_1 + m\omega_2$. Indeed, the most prominent peaks appear at $\omega_1 \pm \omega_2$, followed by $2\omega_1 \pm \omega_2$, confirming the quasiperiodic mechanism.

\FloatBarrier

\section{Master Stability Function Analysis}
\label{sec:msf_analysis}

In this section, we compute the Master Stability Function (MSF) for the system configuration defined in \eqref{Network specific} with parameters $\alpha=2.3427$, $\gamma=0.5$, $D_u=D$, and $D_v=D_w=0$.
As established in Sec.~\ref{Networks with complex eigenvalues}, the stability is determined by the complex vector $\Omega = \lambda^s(D,0,0)^T$. Thus, in the present setting the full MSF parameter space is sampled only along a one-dimensional subspace. For notational simplicity, we denote the nonzero component of this vector by $\Omega=D\lambda^s$. With this convention, the variational equation can be written in terms of the scalar complex parameter $\Omega$. Using the spectral symmetry, we restrict our analytical derivation to $\lambda^s$ in the third quadrant without loss of generality. However, to illustrate the complete stability landscape, the numerical results are presented across both the second and third quadrants, obtained by reflecting the numerical data.

We calculate the Floquet multipliers of the linearized variational system:
\begin{equation}
    \begin{pmatrix} \dot{U}_j \\ \dot{V}_j \\ \dot{W}_j \end{pmatrix} = 
    \begin{pmatrix}
        \Omega + \gamma - 2u_0(t) - \alpha v_0(t) & -\alpha u_0(t) & 0 \\
        0 & 1 - 2v_0(t) - \alpha w_0(t) & -\alpha v_0(t) \\
        -\alpha w_0(t) & 0 & 1 - 2w_0(t) - \alpha u_0(t)
    \end{pmatrix}
    \begin{pmatrix} U_j \\ V_j \\ W_j \end{pmatrix},
    \label{eq:linear_omega}
\end{equation}
and identify the region in the complex $\Omega$-plane where the maximal Floquet multiplier satisfies $|\mu_1| < 1$.

It is convenient to express the complex parameter in polar coordinates relative to the negative real axis:
\begin{equation}
    \Omega = -R e^{i\theta}, \quad R \ge 0, \quad 0 \le \theta < \frac{\pi}{2}.
\end{equation}
Note that under this definition, $\text{Re}(\Omega) = -R\cos\theta$ and $\text{Im}(\Omega) = -R\sin\theta$, which places $\Omega$ in the third quadrant as required.

\begin{figure}[htbp]
    \centering
\includegraphics[width=0.35\linewidth]{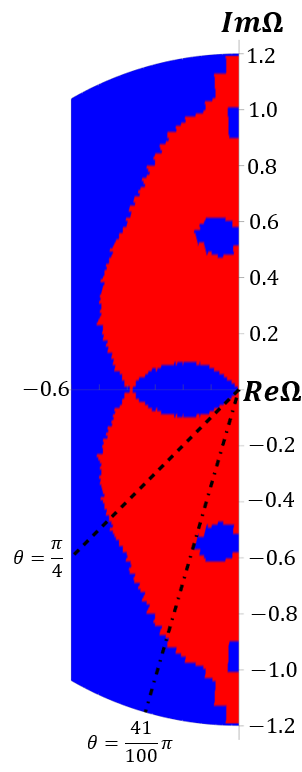}
\caption{Master Stability Function in the complex $\Omega$-plane for the network variational problem, computed numerically with $\alpha=2.3427$, \(\gamma=0.5\), and $D_v=D_w=0$. 
Blue denotes the stable region $(|\mu_{1}|<1)$, and red denotes the unstable region $(|\mu_{1}|>1)$. 
The relevant part of the left half-plane is shown, in accordance with the spectral constraints of zero-row-sum Metzler coupling matrices. 
The dashed rays in the lower half-plane mark characteristic angles measured from the negative real axis: $\theta=\pi/4$, corresponding to the asymptotic transition obtained from the small-$|\Omega|$ expansion, and $\theta=41\pi/100$, corresponding to the numerically determined angle at which the instability region changes its radial structure and multiple instability intervals appear.}
    \label{fig:MSF}
\end{figure}

Numerical analysis reveals the following geometric properties of the instability region (where $|\mu_1| > 1$):
\begin{enumerate}
    \item Fixing the angle $\theta$ and varying the modulus $R$, we find that for sufficiently small $\theta$, the instability region is a single connected interval $I_\theta \subset \mathbb{R}_+$.
    \item The instability intervals exhibit an inclusion property: for $\theta_1 < \theta_2$, we observe $I_{\theta_1} \subset I_{\theta_2}$.
    \item For angles $\theta > \frac{41}{100}\pi$, the topology changes; multiple disjoint instability intervals may appear, and the strict inclusion property no longer holds.
    \item For $R > 1.2$, the system is stable for all angles $\theta$.
\end{enumerate}
The complete MSF landscape is visualized in Figure \ref{fig:MSF}.

To understand the instability mechanism near the origin, we perform an asymptotic expansion of the leading multiplier $\mu_1$ for small $|\Omega|$. Based on the diffusion limit ($\theta=0$), we posit:
\begin{equation}
    \mu_1(\Omega) = 1 + \mu_1^{(1)}\Omega + \mu_1^{(2)}\Omega^2 + o(\Omega^2),
\end{equation}
where numerical evaluation yields $\mu_1^{(1)} \approx 35.7 > 0$ and $\mu_1^{(2)} \approx -5.8 \cdot 10^3 < 0$.
Substituting $\Omega = -R(\cos\theta + i\sin\theta)$, we obtain:
\begin{equation}
    \mu_1(R,\theta) = 1 - \mu_1^{(1)}R\cos\theta + \mu_1^{(2)}R^2\cos 2\theta + i\left( -\mu_1^{(1)}R\sin\theta + \mu_1^{(2)}R^2\sin 2\theta \right) + o(R^2).
\end{equation}
The squared modulus is given by:
\begin{equation}
\label{R*}
    |\mu_1(R,\theta)|^2 = 1 - 2\mu_1^{(1)}R\cos\theta + \left( (\mu_1^{(1)})^2 + 2\mu_1^{(2)}\cos 2\theta \right)R^2 + o(R^2).
\end{equation}
Notably, the ratio of coefficients is small: $\mu_1^{(1)} / |\mu_1^{(2)}| \sim 6.2\cdot 10^{-3}$. This explains why the quadratic term has a visible effect already for moderately small values of \(R\).
Neglecting the term linear in $R$, the stability condition can be approximated as $2\mu_1^{(2)}\cos 2\theta < 0$. Since $\mu_1^{(2)} < 0$, stability requires $\cos 2\theta > 0$, i.e., $\theta < \pi/4$.
This analytical prediction aligns with the numerical results in Figure \ref{fig:MSF}: the system is stable for $\theta < \pi/4$ and unstable for $\theta > \pi/4$ (in the intermediate $R$ regime).
However, from \eqref{R*} it can be seen that when $\theta > \pi/4$, the solution remains stable within a vanishingly small interval $0 \le R \le R_*$, where:
\begin{equation}
    R_* \approx \frac{2\mu_1^{(1)}\cos\theta}{2\mu_1^{(2)}\cos 2\theta + (\mu_1^{(1)})^2}.
\end{equation}
This interval is practically invisible in standard plots but ensures consistency with the stability of the uncoupled limit ($R=0$).

\FloatBarrier 

\newpage

\section{Conclusion}

We have shown that finite-wavenumber instabilities in continuous media can be
mapped into spectral instability windows in network coupling space.
We established a linear spectral correspondence between continuous LSI media and discrete networks, showing that Fourier modes of a continuous spatial operator and transverse eigenmodes of a network coupling matrix enter the corresponding variational equations in the same spectral form. 
In this sense, the Master Stability Function provides a common stability plane for both continuous and discrete systems.

In the present work, this correspondence was used primarily in a predictive direction. 
Starting from a finite-wavenumber instability of a continuous
reaction-diffusion extension, the eigenvalues of the network coupling matrix
determine the coupling intervals in which synchronization is lost and recovered. 
This yields a mechanism for reentrant synchronization: increasing coupling strength does not necessarily stabilize the synchronized state monotonically, but may instead drive the network through alternating synchronization and desynchronization windows.

Using the competitive three-species Lotka-Volterra model as a representative example, we demonstrated this mechanism for networks with zero-row-sum Metzler coupling matrices. 
For real transverse eigenvalues, the instability mechanism is inherited directly from the finite-wavenumber period-doubling instability of the continuous problem. 
For complex transverse eigenvalues, which arise naturally in directed networks, the instability threshold is crossed by complex Floquet multipliers. 
The resulting bifurcation is quasiperiodic rather than period-doubling, producing a dynamical regime that is absent in standard reaction-diffusion systems with real wavenumbers.

The analysis also shows the importance of the algebraic structure of the coupling matrix. 
Zero-row-sum Metzler matrices, which naturally describe conservative transport and population exchange between habitats, restrict the relevant spectrum to sectors of the left complex half-plane. 
These constraints reduce the region of the MSF plane that must be examined and provide a useful link between network topology and synchronization stability.

The correspondence developed here is a linear spectral-stability correspondence and does not imply a nonlinear equivalence between continuous and discrete systems. 
Nevertheless, it provides a systematic way to transfer intuition from continuous dispersion relations to network spectra, and to identify where this analogy breaks down. 
In particular, directed networks can access complex spectral directions that are not available in standard real-wavenumber reaction-diffusion models.

Several questions remain open. 
The present work focused on the onset of linear instability and on isolated instability windows. 
The nonlinear evolution of the unstable modes, the interaction between overlapping instability intervals, and the coexistence of period-doubling and quasiperiodic mechanisms require further investigation. 
Numerical continuation and direct simulations of the full nonlinear network dynamics may provide insight into these regimes.

Another natural direction is to extend the correspondence to continuous systems with genuinely complex Fourier symbols, such as advection-diffusion systems or models with asymmetric nonlocal kernels \cite{Klausmeier1999,Siebert2014}. 
Such extensions would clarify the relation between continuous spectral curves and the more flexible spectral sampling available in directed networks. 
It would also help determine which instability mechanisms are genuinely network-induced and which have continuous analogs in more general LSI media. 
Extensions to higher-dimensional spatial domains and to broader classes of coupling operators are also promising directions for future work.

\section{Acknowledgments}
The authors thank Alvin Bayliss and Michael Zaks for fruitful discussions and valuable suggestions.

\appendix

\section{The non-diagonalizable case}
\label{Appendix A}

Consider the case where the coupling matrix $C$ is defective (not diagonalizable) and is represented by its Jordan normal form. In this scenario, the perturbations cannot be expanded solely in terms of eigenvectors. Instead, we employ a basis of generalized eigenvectors.

Let $\{\mathbf{v}^{(1)}, \dots, \mathbf{v}^{(l)}\}$ be a Jordan chain of length $l$ associated with the eigenvalue $\lambda$, satisfying the relations $(C - \lambda I)\mathbf{v}^{(1)} = 0$ and $(C - \lambda I)\mathbf{v}^{(j)} = \mathbf{v}^{(j-1)}$ for $j=2,\dots,l$.
Expanding the perturbation $\delta \mathbf{u}(t)$ in this basis, the variational equations for the mode amplitudes $\boldsymbol{\xi}_j(t)$ become coupled in a cascade structure:

\begin{subequations}
\label{eq:jordan_full}
\begin{align}
\frac{d\boldsymbol{\xi}_1}{dt} &= \left[ D\mathbf{f}(\mathbf{u}_0) + \lambda \mathbf{\tilde{D}} \right] \boldsymbol{\xi}_1, \label{eq:jordan_base} \\
\frac{d\boldsymbol{\xi}_j}{dt} &= \left[ D\mathbf{f}(\mathbf{u}_0) + \lambda \mathbf{\tilde{D}} \right] \boldsymbol{\xi}_j + \mathbf{\tilde{D}} \boldsymbol{\xi}_{j-1}, \quad \text{for } j=2,\dots,l. \label{eq:jordan_chain}
\end{align}
\end{subequations}

Crucially, the equation for the first vector in the chain \eqref{eq:jordan_base} is identical to the standard Master Stability Function equation \eqref{eq:network_variational}.
This implies that the non-diagonalizable nature of the network does not eliminate instabilities. If the continuous LSI system is unstable for a spectral parameter $\hat{g}(\mathbf{k}) = \lambda$, then the network dynamics will also be unstable. Specifically, the mode $\boldsymbol{\xi}_1$ will grow exponentially at the same rate as in the diagonalizable case. Furthermore, the subsequent modes $\boldsymbol{\xi}_j$ will be driven by the instability of $\boldsymbol{\xi}_{j-1}$, potentially leading to algebraic growth factor ($t^{j-1} e^{\text{Re}(\sigma)t}$) on top of the exponential growth. Thus, the instability criterion derived from the continuous spectrum remains a sufficient condition for instability in general network topologies.
If the MSF is strictly stable at the eigenvalue, Jordan blocks may introduce
only polynomial prefactors and therefore do not destroy exponential stability.
However, at marginal points, where the leading Floquet exponent has zero real
part or a multiplier lies on the unit circle, the Jordan structure may affect
the local bifurcation and requires separate analysis.

\section{{Attractors of the local model}}
\label{Appendix B}
Consider the Lotka-Volterra system:
\begin{eqnarray}
\label{eq19}
u_t &= u(\gamma - u - \alpha v) \nonumber\\ 
v_t &= v(1 - v - \alpha w)\\ 
w_t &= w(1 - w - \alpha u)\nonumber  
\end{eqnarray}
with the growth rate of the first species $\gamma \neq 1$.   

\subsection{Linear instability of the coexistence state}

The coexistence solution is
\begin{equation}
\label{eq20}
u_*=\frac{\gamma -\alpha +\alpha ^2}{1+\alpha ^3},\; 
v_*=\frac{1 -\alpha +\alpha ^2 \gamma}{1+\alpha ^3},\; 
w_*=\frac{1 -\alpha \gamma +\alpha ^2}{1+\alpha ^3}.
\end{equation}

We assume that
\begin{equation}
\label{eq21} 
\alpha >1,\; \frac{\alpha -1}{\alpha ^2}<\gamma <\frac{\alpha ^2 +1}{\alpha},
\end{equation}
hence $u_*$, $v_*$, $w_*$ are positive. 

Stability of the coexistence state is determined by the linearization of (\ref{eq19});
three eigenvalues $\sigma$ of the Jacobian matrix solve the equation
\begin{equation}
\label{eq23}
\sigma ^3+(u_*+v_*+w_*)\sigma ^2+(u_*v_*+v_*w_*+w_*u_*)\sigma+(1+\alpha ^3)u_*v_*w_*=0
\end{equation}

The stability boundary is when $\sigma=i\omega$ for a real $\omega$. Substituting this form of $\sigma$ into (\ref{eq23}) gives two equations, on the real and imaginary parts:
$-\omega ^2 +u_*v_*+v_*w_*+w_*u_*=0 \ , \\ -\omega ^2(u_*+v_*+w_*)+(1+\alpha ^3)u_*v_*w_*=0$, therefore:

\begin{equation}
\label{eq47}
-(u_*v_*+v_*w_*+w_*u_*)(u_*+v_*+w_*)+(1+\alpha ^3)u_*v_*w_*=0.
\end{equation}

Substituting \eqref{eq20} into \eqref{eq47} determines a curve on the plane $(\alpha ,\gamma)$:

\begin{equation}
\label{eq25}\begin{split}
\alpha (1 - \alpha - \alpha^3) \gamma^3 
+ (-2 + 3\alpha - 5\alpha^2 + 6\alpha^3 + \alpha^5 - \alpha^6) \gamma^2\\ 
+ (-4 + 7\alpha - 11\alpha^2 + 5\alpha^3 - 7\alpha^4 + \alpha^5 + \alpha^7) 
\gamma\\ 
+ (1 - \alpha)^2 (-2 + \alpha - 3\alpha^2 - \alpha^3 - \alpha^4) = 0.
\end{split}
\end{equation}
The curve of Hopf bifurcation splits the parameter plane into the region 
where the coexistence point is an attractor, and the region where it is unstable.

\subsection{Stability of the heteroclinic cycle}

System (\ref{eq19}) has three saddle points $(\gamma,0,0)$, $(0,1,0)$, $(0,0,1)$, which are connected by heteroclinic trajectories forming a heteroclinic cycle.

Consider a trajectory that is very close to the heteroclinic cycle.
The time needed to complete a "cycle" (start close to one of the saddle points and then come back to its neighborhood) can be divided into 6 intervals: 3 intervals where the solution is close to one of the saddles, which means that only one of the variables is $O(1)$, and 3 intervals on each of the heteroclinic trajectories, where two of the variables are $O(1)$. We show that the major contribution comes from the time the system spends near the saddles.

Consider the trajectory starting between the vicinities of saddle points $(0,0,1)$ and $(\gamma,0,0)$, where $u$ and $w$ are $O(1)$, while $v=\epsilon \ll 1$.
Near the saddle point $(\gamma,0,0)$, the system is governed by the linear equations:
$$\tilde{u}_t=-\gamma \tilde{u}-\alpha \gamma v, v_t=v, w_t=(1-\alpha \gamma )w.$$
Here $\tilde{u}=u-\gamma$.
The time $t_1$ spent near that point depends on $\epsilon$: during that time, the variable $v$ increases from $\epsilon$ to $O(1)$. Since $\Tilde{v} \sim \epsilon e^t$, we have $\epsilon e^{t_1} \sim 1$, so $t_1 \sim -\ln\epsilon$. 
During that time, $w$ decreases from $O(1)$ to $\Tilde{w_1} \sim e^{(1-\alpha \gamma)t_1}$. 

When $v$ becomes $O(1)$, the linear equation for $v$ has to be corrected.
Between the saddles $(\gamma,0,0)$ and $(0,1,0)$, where both $v$ and $1-v$ are $O(1)$, the governing equation for $v$ is $v_t=v(1-v)$ and the solution for $v$ corresponds to the heteroclinic trajectory in the plane $w=0$,
\begin{equation}
\label{eq26}
v(t)=\frac{1}{1+C_1e^{-t}},
\end{equation}
where $C_1$ is a constant depending on the initial conditions. Because for $v=O(1)$ and $1-v=O(1)$, the velocity $v_t=O(1)$, the duration of the motion on that fragment of trajectory is $O(1)$. Hence, it is much smaller than $t_1$, and the variable $w$ does not change its order during that interval, hence the change of the variable $w$ during that time can be neglected.
Near the point $(0,1,0)$, the evolution of variables is governed by the system of equations $$u_t=(\gamma-\alpha)u, \ \tilde{v_t}=-\tilde{v}-\alpha w,\; w_t=w,$$
where $\tilde{v}=v-1$.
During the time $t_2$ spent near that point, $w$ grows from $\Tilde{w_1}$ up to $O(1)$, hence, similarly to the calculation of $t_1$ we get: $t_2\sim -\ln(\Tilde{w_1})=-(\alpha \gamma -1)\ln\epsilon$.
During that time $t_2$, $u$ decreases from $O(1)$ to $u_2\sim e^{(\gamma-\alpha)t_2}=\epsilon^{(\alpha-\gamma)(\alpha \gamma -1)}$.

Between the points $(0,1,0)$ and $(0,0,1)$ the governing equation for $w$ is $w_t=w(1-w)$,
and the heteroclinic solution on the plane $u=0$  is
\begin{equation}
\label{eq27}
w(t)=\frac{1}{1+C_2e^{-t}},
\end{equation}
where $C_2$ is a constant depending on the initial conditions. Again, the duration of the time interval where $w=O(1)$ and $1-w=O(1)$ is $O(1)$, and the variable $u$ does not change its order during that interval.

Near the point $(0,0,1)$, the evolution of variables is governed by the system of equations,
$$u_t=\gamma u, \ v_t=(1-\alpha)v, \ \tilde{w_t}=-\tilde{w}-\alpha \tilde{u},$$
where $\tilde{w}=w-1$.
The variable $u$ grows from 
$\Tilde{u_2}$ to $O(1)$ during the time $t_3$: $t_3\sim -\frac{1}{\gamma}\ln(\Tilde{u_2})=-\frac{(\alpha -\gamma)(\alpha \gamma -1)}{\gamma}\ln\epsilon$. \\
During that time $v$ decreases from $O(1)$ to $\Tilde{v_3}\sim e^{-(\alpha-1)t_3}=\epsilon^{(\alpha-1)(\alpha-\gamma)(\alpha \gamma-1)/\gamma}$. \\
Finally, between the points $(0,0,1)$ and $(\gamma,0,0)$ the governing equation for $u$ is $u_t=u(\gamma-u)$ and the solution for $u$ corresponds to the heteroclinic trajectory in the plane $v=0$,
\begin{equation}
\label{eq28}
u(t)=\frac{\gamma}{1-e^{\gamma(C_3-t)}}
\end{equation}
where $C_3$ is a constant depending on the initial conditions.
The variable $v$ does not change its order in that interval.

We find that during the motion in the vicinity of the heteroclinic cycle $(\gamma,0,0)\rightarrow (0,1,0)\rightarrow (0,0,1)\rightarrow (\gamma,0,0)$ the order of the variable $v$ is changed from $\epsilon$ to $\epsilon^c$ where 
$$c=(\alpha-1)(\alpha-\gamma)(\alpha \gamma -1)/\gamma.$$ 
When $c>1$, the heteroclinic cycle is an attractor since a solution that starts close to the heteroclinic cycle becomes even closer to it after one cycle. When $c<1$, the heteroclinic cycle is a repeller. Note that if we start from a different saddle point, we obtain the same coefficient $c$ (with the same factors appearing in a different order). The critical case $c=1$ divides the plane $(\gamma,\alpha)$ by the curve:
\begin{equation}
\label{eq29}
(\alpha-1)(\alpha \gamma -1)(\alpha -\gamma)-\gamma=0.
\end{equation}
Equations (\ref{eq25}),(\ref{eq29}) together with inequalities (\ref{eq21}) split the relevant
parameter plane $(\gamma,\alpha)$ into three regions: one where the fixed point is an attractor, one where the heteroclinic cycle is an attractor, and an
intermediate region between them. This suggests the existence of an attracting limit cycle in the intermediate
region. This is consistent with the analysis of the attractors of a three-species competitive Lotka-Volterra model \cite{Coste1979}.

\end{document}